\def\ETDS{1}
\if0=1
\documentclass{etds}
\newtheorem{theo}{Theorem} 
\newtheorem{lemma}{Lemma}

\begin{document}
\ETDS{1}{5}{XX}{1999}

\runningheads{A.\ Prikhod'ko}{Lower bounds for symbolic complexity of iceberg dynamical systems}

\title{Lower bounds for symbolic complexity of iceberg dynamical systems}

\author{Alexander Prikhod'ko}

\address{Dept.\ of Mechanics and Mathematics, Moscow State University, Moscow, 
Russia\\
\email{sasha.prihodko@gmail.com}}

\recd{$27$ January $2012$}
\else
\documentclass{amsart}
\usepackage{amssymb}
\usepackage{graphics}
\textwidth=164mm
\textheight=215mm
\topskip=0mm
\topmargin=0mm
\oddsidemargin=0mm
\evensidemargin=0mm
\tolerance=8400
%
\newtheorem{theo}{Theorem} 
\newtheorem{lemma}{Lemma}
\newtheorem{prop}{Proposition}
\newtheorem{cor}{Corollary}
\def\proc#1{\addvspace{6pt}{\it #1\/}$\;$ }
\def\ep{\noindent{\hfill $\Box$}}
\def\acks{\proc{Acknowledgements.}} 
\title{Lower bounds for symbolic complexity of iceberg dynamical~systems}
\author{Alexander Prikhod'ko}
\address{Dept.\ of Mechanics and Mathematics, Moscow State University}
\email{sasha.prihodko@gmail.com}
\begin{document}
\fi

%
\def\Set#1{{\mathbb{#1}}}
\def\const{\mathop{\mathrm{const}}\nolimits}
\def\supp{\mathop{\mathrm{supp}}\nolimits}
\def\Mix{\mathop{\mathrm{Mix}}\nolimits}
\def\WMix{\mathop{\mathrm{WMix}}\nolimits}
\def\Id{\mathop{\mathrm{Id}}\nolimits}
\def\id{\mathop{\mathrm{id}}\nolimits}
\def\eset{\varnothing}
\def\Maps{\colon\:}
\def\where{\colon\:}
\def\impl{\;\Rightarrow\;}
\def\wt#1{\widetilde{#1}}
\def\wto{\mathrel{\stackrel{\mathrm w}{\to}}}
\def\sms{\mathbin{\smallsetminus}}
\def\syms{\mathbin{\triangle}}
\def\df{{\boldsymbol\delta}}
\def\dbar{{\Bar d}\,}
\def\No{No.\,}
\def\a{\alpha}
\def\eps{\varepsilon}
\def\la{\lambda}
\def\cB{{\mathcal{B}}}
\def\cD{{\mathcal{D}}}
\def\cE{{\mathcal{E}}}
\def\cL{{\mathcal{L}}}
\def\cP{{\mathcal{P}}}
\def\cV{{\mathcal{V}}}
%
\def\om{\alpha}
\def\p{p}
\def\pt{{\boldsymbol \pi}}
\def\ident{\mathrel{\stackrel{\bullet}=}}
\def\ptjoined{\mathrel{\leftrightarrow}}
\def\sC{\bar\p}
\def\ito{\,..\,}
\def\iw{{\boldsymbol w}}
\def\tW{{\Tilde W}}
\def\Lang{\cL}
\def\cyclsub{\mathbin{\preceq_\circ}}
\def\Prob{{\mathsf P}}
\def\scale{{\boldsymbol\la}}

\begin{abstract}
The symbolic complexity of an infinite word $W$ is the function $\p_W(l)$ 
counting the number of different subwords in~$W$ of length~$l \in \Set{N}$. 
In this paper our main purpose is to study the complexity for a class of 
topological dynamical systems, called {\it iceberg systems}, 
given by the following symbolic procedure. 
Starting from a given finite word $\iw_1$ 
we construct a sequence of words 
$$\iw_{n+1} = \iw_n\, \rho_{\om_n(1)}(\iw_n)\cdots\rho_{\om_n(q_n-1)}(\iw_n),$$ 
where $\rho_\a(\iw_n)$ is a~cyclic rotations of the word~$\iw_n$, and 
consider an infinite word $\iw_\infty$ extending each $\iw_n$ to the right. 
It is shown that for iceberg systems given by the randomized parameters $\om_n(j)$ 
the complexity function satisfies the estimate 
${ \p_{\iw_\infty}(l) \gtrsim l^{3-\eps} }$ for any ${\eps > 0}$, 
and at the same time it is observed that this estimate represents up to a small correction the optimal 
lower bound for the complexity function, namely, 
${ \p_{\iw_{n+1}}(l_n) \le l_n^3 }$ along the subsequence ${l_n = |\iw_n|+1}$.  

The work is supported by grant NSh-8508.2010.1 for support of leading Russian scientific schools  
and RFFI grant \No\,11-01-00759-a. 
\end{abstract}

\maketitle

\section{Introduction}

The symbolic complexity of an infinite word $W$ is the function $\p_W(l)$ 
counting the total number of different subwords in~$W$ of length~$l \in \Set{N}$. 
In a~series of recent investigations the concept of symbolic complexity 
is applied to the study of orbit structure both in topological dynamics 
and ergodic theory of measure preserving transformations 
(e.g.\ see \cite{ArnouxMauduitShTam,BRL,delVAR,Ferenczi2,FerencziZamboniOnIET08,MelaPetersenOnDynPrPascal}). 
Further, a~kind of invariant of measure-theoretic isomorphism 
extending the concept of symbolic complexity 
is introduced by S.\,Ferenczi in \cite{FerencziOnMTCompl}. 
Another general approach to the complexity or orbit structure 
introduced by A.\,Vershik is based on the notion 
of scaling entropy \cite{VershikOnScalingEntrAndPPS}.

Out investigation is motivated by the idea of applying the concept of complexity 
and coding arguments to the problems of isomorphism in ergodic theory of 
measure preserving transformations with zero entropy. 
We would like to start with mentioning the classical symbolic construction 
of a {\it rank one\/} dynamical system. 
Given a non-trivial finite word ${\boldsymbol v}_{1}$ in a finite alphabet $\Set{A}$, 
we define words 
\begin{equation}\label{eRankOneSymb}
	{\boldsymbol v}_{n+1} = 
		{\boldsymbol v}_n \: 0^{s_{n,0}} \: {\boldsymbol v}_n \: 0^{s_{n,1}} \: 
		{\boldsymbol v}_n \: \cdots \: 0^{s_{n,r_n}} \: {\boldsymbol v}_n,   
\end{equation}
where symbol ``$0$'' used to create small {\it spacers\/} between the copies of the word~${\boldsymbol v}_n$, 
and an infinite word ${\boldsymbol v}_\infty$ having each ${\boldsymbol v}_{n}$ as a prefix. 
Following the standard scheme we can consider the compact set $K$ defined as the weak closure 
of all translates of ${\boldsymbol v}_\infty$ and endow $K$ with the standard measure 
which is invariant under the shift trandformation $T$. 
The triple $(T,K,\mu)$ is called {\it rank one\/} transformation. 

Rank one transformation serves as a classical example of dynamical system with simple spectrum. 
In~paper \cite{IcePaperI} a class of ergodic maps is introduced, 
which is similar to rank one systems and extends this class. 
Namely, generalizing the rule~\eqref{eRankOneSymb}, 
while defining $\iw_{n+1}$ we combinate not only exact copies but 
any cyclic rotations of the previous word $\iw_n$. 
It is observed that the measure preserving map associated with this new symbolic procedure, 
called {\it iceberg transformation},  
has always $1/4$-local rank and, under certain conditions, simple spectrum. 
It is not obvious that the new class is actually larger than the class of rank one transformations. 
At the same time, these two classes potentially could contain transformations 
having spectral types with similar properties. 
Thus, it is interesting to find an invariant helping to distinguish these classes, 
and, in fact, it can be easily deduced from the results of this paper that 
the typical complexity of iceberg systems is $\p_{\iw_\infty}(l) \gtrsim l^{3-\eps}$, hence, 
a wide class of iceberg systems is certainly goes beyond the class of rank one systems, 
since for any rank one map there exists a finite partition of the phase space 
generating a word ${\boldsymbol v}_\infty$ with the property 
$\p_{{\boldsymbol v}_\infty}(l_n) \lesssim \frac12 l_n^2$ (see~\cite{Ferenczi1}) 
in contradiction with the given lower bound, 
proving the existence of iceberg systems which are not rank one. 
Though, one can observe that this approach is essentially non-spectral, hence, 
there are no  obstacles to existence of a pair $(T_1,T_2)$ of spectrally isomorphic maps, 
where $T_1$ is rank one and $T_2$ is an iceberg map, but not rank one.
Thus, this class of symbolic systems, more general than rank one, can be considered 
as a~source of new examples demonstrating that rank is not a spectral invariant. 

It~is interesting to compare the lower boud for symbolic complexity of iceberg systems 
with the asymptotics of the Pascal adic transformation complexity, ${p(l) \sim l^3/6}$ 
established by X.\,M\'ela and K.\,Petersen \cite{MelaPetersenOnDynPrPascal}.


\section[Iceberg systems]{Iceberg dynamical systems}

We start with the a series of constructions and definitions. 

\proc{Definition.}
Consider a finite alphabet $\Set{A}$ and a word ${w = a_1\dots a_N}$ in~$\Set{A}$. 
We use the notation $|w|$ for the length of the wor~$w$. 
Given an integer number~$\a$, ${0 < \a < |w|}$, let us define the 
{\it cyclic rotation\/} $\rho_\a(w)$ of the word~$w$ as follows   
\begin{equation*}
	\rho_\a(w) = w_2 w_1, \quad \text{whenever} \quad w = w_1 w_2, \quad |w_1| = \a. 
\end{equation*}
%
Observe that if $b \in \Set{A}$ is a~letter and $u$ is a~word 
then ${\rho_1(bu) = ub}$ and ${\rho_\a = (\rho_1)^\a}$. 
Let us 
extend this definition to all values of ${\a \in \Set{Z}}$ setting 
\begin{equation*}
	\rho_0(w) = w \quad \text{and} \quad \rho_\a(w) = \rho_{\a+k|w|}(w), \quad k \in \Set{Z}. 
\end{equation*}
\medbreak

The next definition is the starting point of the construction of iceberg systems. 
We start with defining an infinite word $\iw_\infty$, then we consider the closure $K_{\iw_\infty}$ 
of the sequence $\{\iw_\infty\}$ in weak topology, and finally using the sandard procedure 
we endow the compact set $K_{\iw_\infty}$ with an invariant measure $\mu^{(\iw_\infty)}$ 
and we come to an ergodic transformation $T$ acting as left shift 
on the space $(K_{\iw_\infty},\mu^{(\iw_\infty)})$ (see~\cite{IcePaperI}). 

\proc{Definition.}\label{defIcebergTopSystWOS}
Let $\iw_1$ be a fixed finite word in~$\Set{A}$. 
We suppose that $\iw_1$ is non-trivial and contains a pair of different letters. 
Let us construct 
a sequence of words $\iw_n$, such that each next word $\iw_{n+1}$ in the sequence 
is a concatenation of cyclic rotations of $\iw_n$, namely, 
\begin{equation}
	\iw_{n+1} = \rho_{\om_n(0)}(\iw_n)\, \rho_{\om_n(1)}(\iw_n) \cdots \rho_{\om_n(q_n-1)}(\iw_n), 
\end{equation}
where $q_n \ge 2$ is the number of entries $\rho_{\om_n(j)}(\iw_n)$. 
Now suppose, for simplicity, 
that ${\om_n(0) = 0}$ for any~$n$ so that $\iw_{n+1}$ extends~$\iw_n$ 
to the right, and after the infinite sequence of steps we come to 
an infinite word~$\iw_\infty$ such that any $\iw_n$ is prefix of~$\iw_\infty$. 
Let $K_{\iw_\infty}$ be the closure of $\{\iw_\infty\}$ in weak topology, 
i.e.\ the set of all sequences $x = (x_j)_{j\in \Set{Z}}$ such that any finite subword of~$x$ 
can be found as a subword of some $\iw_n$. 
The compact set $K_{\iw_\infty}$ is invariant under the left shift map 
\begin{equation}
	T \Maps (\dots,x_0,x_1,\dots,x_j,\dots) \mapsto (\dots,x_1,x_2,\dots,x_{j+1},\dots). 
\end{equation} 
We~call $(T,K_{\iw_\infty})$ an {\it iceberg topological dynamical system\/} 
({\it without spacers\/}). 
\medbreak


\subsection{Symbolic complexity.}

Given an infinite word $w$ let us consider the {\it language\/} $\Lang(w)$ 
defined as the set of all finite subwords of the word~$w$. 
We can partition the set $\Lang(w)$ 
into a sequence of subsets $\Lang(w,l)$ according to the length of a subword, 
\begin{equation*}
	\Lang(w,l) = \{u \in \Lang(w) \where |u| = l\}. 
\end{equation*}
Symbolic complexity $\p_w(l)$ of the word~$w$ is the function 
counting the number of different subwords in~$w$ of length~$l$ or, in other term, 
the function measuring the complexity of $\Lang(w)$, 
\begin{equation*}
	\p_w(l) = \#\Lang(w,l), 
\end{equation*}
where $\# A$ is the cardinality of~$A$. 

Our paper is devoted to the proof of the following observation. 

\begin{theo}
Consider an iceberg system given by the independent uniformely distributed parameters $\a_n(j)$, 
and assume that ${q_n \to \infty}$ sufficiently fast, ${q_n \ge h_n^\gamma}$ with ${\gamma \ge 1}$, 
where ${h_n = |\iw_n|}$. 
Then almost surely 
\begin{equation}
	\p_{\iw_\infty}(l) \ge l^{3-\eps} 
\end{equation}
for any ${\eps > 0}$ and $l \ge l_0(\eps)$. 
\end{theo}


\subsection{Measure-theoretic point of view.}

Throughout this paper we will continuously compare pure combinatorial observations with the 
corresponding effects in the context of measure-theoretic ergodic theory. 
Infinite words are interpreted in this case as discrete sample path of stationary random processes. 

\begin{figure}[th]
  \centering
  \unitlength=1mm
  \includegraphics{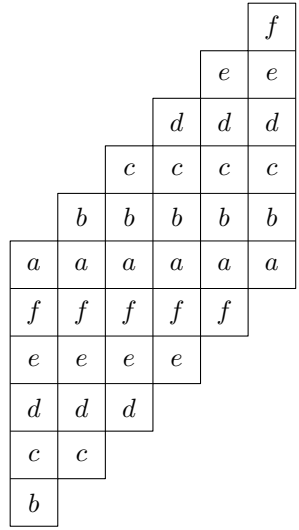} 
  \caption{Sample of an iceberg. 
  	A~column is an ordinary Rokhlin tower corresponding to a~cyclic rotation of the word ``$abcdef$''. 
  	Transformation $T$ lifts each set drawn as a square to the upper square. 
  	Whenever a point reaches the top set in a~column, it is mapped by $T$ to an arbitrary bottom set 
  	of another (or~the~same) column.} 
  \label{fJumpsOnXn}
\end{figure}

The concept of iceberg transformation can be extended 
if we endow the compact space $K_{\iw_\infty}$ 
with the natural invariant measure $\mu = \mu^{(\iw_\infty)}$ 
generated by the set of empirical distributions $\mu_{l}$ 
definied as~follows: given a~finite word ${w \in \Lang(\iw_\infty,l)}$, a~word of length~$l$, 
the value $\mu_{l}(w)$ is the asymptotic frequency of occurence of~$w$ inside the infinite word~$\iw_\infty$. 
This procedure leads to an ergodic dynamical system given by 
the left shift map $T$ acting on the probability space $(K_{\iw_\infty},\cB,\mu^{{(\iw_\infty)}})$, 
where $\cB$ is the sigma-algebra of Borel sets. 
It is shown in  paper~\cite{IcePaperI} that an iceberg map can be obtained as a~result 
of cutting-and-stacking construction in a~way independent on the symbolic formalism we discuss 
(see~fig.\,\ref{fJumpsOnXn}). 

The main idea of our paper is to provide a combinatorial background to 
further study of a~typical symbolic complexity of ergodic iceberg systems 
with respect to the distribution~$\mu^{(\iw_\infty)}$ 
(cf.~\cite{FerencziOnMTCompl} and \cite{VershikOnScalingEntrAndPPS}). 

When studying our dynamical system as an ergodic transformation it is naturally to ask 
if the underlying property of the iceberg construction can be formulated in an invariant form, 
as a kind of dynamical property preserved under any measure-theoretic isomorphism? 

Indeed, let us consider the classical concept of rank one approximation. 

\proc{Definition.}
We say that a measure preserving map $T$ on a probability space $(X,\mu)$ 
is of {\it rank one\/} if for any measurable partition ${\cP = \{P_1,\dots,P_m\}}$ 
and any ${\eps > 0}$ there exists a set $\Omega_\eps$ 
of $\mu$-measure ${1-\eps}$ such that for any ${x \in \Omega_\eps}$ the orbit of $x$ 
represented in a form of an infinite word $(x_j)$, where ${T^j x \in P_{x_j}}$, 
satisfies the following property: $(x_j)$ is $\eps$-covered by disjoint words $\tW_j$ 
which are $\eps$-close to some fixed word $W_{(\eps)}$: 
\begin{equation}
	(x_j) = \cdots \tW_1 \,S_1\, \tW_2 \,S_2\, \cdots \,S_{j-1}\, \tW_j \cdots, 
\end{equation}
where $\tW_j$ concides with $W_{(\eps)}$ up to at most $\eps|W_{(\eps)}|$ wrong letters. 
\medbreak

\proc{Definition.}
Let $T$ be a measure preserving transformation. 
We say that $T$ admits {\it iceberg approximation\/} 
if for~any finite measurable partition $\cP$ and any ${\eps > 0}$ one can find a~set $\Omega_\eps$ 
of~measure ${1-\eps}$ such that for any point ${x \in \Omega_\eps}$ the orbit $(x_j)$ of~$x$ 
is $\eps$-covered by words $\tW_j$ which are $\eps$-close to cyclic rotations  
of a~fixed $W_{(\eps)}$: 
\begin{equation}
	(x_j) = \dots \tW_1 \,S_1\, \tW_2 \,S_2\, \dots \,S_{j-1}\, \tW_j \dots, \qquad 
	\dbar\bigl( \tW_j, \rho_{\varphi_j}(W_{(\eps)}) \bigr) < \eps, 
\end{equation}
where $\dbar(u,v)$ measures the fraction of non-matching letters in $u$ and~$v$. 
\medbreak

It is important to say that a~priory we cannot skip spacers in this definition. 
Thus, generally speaking a~common construction of an iceberg map should contain spacers as~well, 
but, at the same time, it is interesting to observe that from the complexity point of view   
adding spacers between cyclic rotations we do not influence significatly to the asymptotics 
of $\p_{\iw_\infty}$ (see the effect discussed below in~remark~\ref{MeltingWordEffectRemark}). 

\section{Upper bound for symbolic complexity along a subsequence}

We start with the estimation of the symbolic complexity 
for the sequence of special kind ${l_n = h_n+1 = |\iw_n|+1}$. 
The idea of this calculation is 
to show that the lower bound discussed in the next section is, in a sense, optimal. 

\begin{lemma}\label{lemUpperBoundIcebergSC}
$\p_{\iw_{n+1}}(h_n+1) \le h_n^3$ for any iceberg system. 
\end{lemma}

\begin{figure}[th]
  \centering
  \unitlength=1mm
  \includegraphics{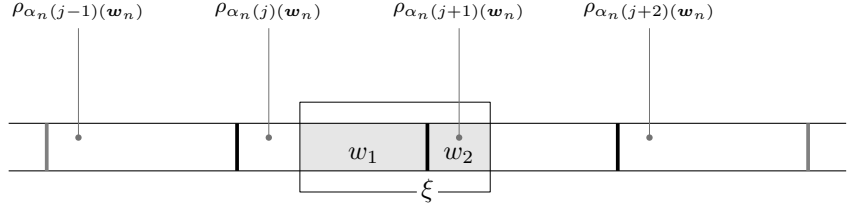} 
  \caption{Estimation of $\p_{\iw_{n+1}}(l)$ for $l$ with ${l = |\iw_n|}$. 
  	The rectangle restricts a subword $v$ of length~$l$. 
  	The vertical lines mark up the boundaries of blocks $\rho_{\om_n(j)}(\iw_n)$.} 
  \label{fUpperBoundIcebergSC}
\end{figure}

\proc{Proof.}
Let us consider a subwords $v$ of length ${h_n+1 = |\iw_n|+1}$ inside $\iw_{n+1}$, 
and consider $\iw_{n+1}$ as a~sequence of rotations $\rho_{\om_n(j)}(\iw_n)$. 
Suppose that $v$ starts at position $k$ so that $v$ includes the letters at positions $k \ito k+h_n$ 
and it covers exactly two adjacent subwords: 
$\rho_{\om_n(j)}(\iw_n)$ and $\rho_{\om_n(j+1)}(\iw_n)$. 
It is easy to see that 
the subword $v$ is uniquely determined by the value of the triple $(\bar k,\om_n(j),\om_n(j))$, 
where ${\bar k = k \pmod{h_n+1}}$, ${\bar k \not= 0}$, 
and the total number of such triples is estimated as~$h_n^3$. 
This ides is illustrated on fig.~\ref{fUpperBoundIcebergSC}. 
Let us summarize the estimates in the following line: 
\begin{equation*}
	\p_{\iw_{n+1}}(h_n) \le h_n^3 - h_n^2 + h_n, \qquad 
	\p_{\iw_{n+1}}(h_n+1) \le h_n^3. \quad \text{\ep} 
\end{equation*}
\medbreak

Now let us consider an iceberg system {\it with spacers\/} given by the scheme 
\begin{equation*}
	\iw_{n+1} = 
		\iw_n \,0^{s_{n,0}}\, 
		\rho_{\om_n(1)}(\iw_n) \,0^{s_{n,1}}\, 
		\cdots \,0^{s_{n,q_n-2}}\, 
		\rho_{\om_n(q_n-1)}(\iw_n) \,0^{s_{n,q_n-1}}, 
\end{equation*}
and observe that spacers 
do~not affect to the asymptotics along the sequence $h_n+1$. 
%
%
We write ${w_1 \preceq w}$ if $w_1$ is a~subword of~$w$. 

\proc{Remark}(Melting word effect). 
\label{MeltingWordEffectRemark}
Consider a subword 
\begin{equation*}
	v = w_1 S w_2, \qquad 
	w_1 \preceq \rho_{\varphi_1}(\iw_n), \qquad 
	w_2 \preceq \rho_{\varphi_2}(\iw_n), 
\end{equation*}
of length $l_n = |\iw_n|+1$ that 
covers a~part $w_1$ of the~block $\rho_{\om_n(j)}(\iw_n)$, 
a~part $w_2$ of the~block $\rho_{\om_n(j+1)}(\iw_n)$, 
and a~spacer ${S = 0^s}$. Let us imagine that the spacer $S$ is growing to the left 
and becomes $\nu$ symbols larger. It can be observed that the growth of the spacer 
can be compensated by the opposite cyclic rotation 
in the first block $\rho_{\om_n(j)}(\iw_n)$ by $-\nu$ positions. 
In other terms, the original word $v$ coincides up to $\nu$~letters with a~new word $\wt v$ 
given by the configuration 
\begin{equation*}
	\wt v = \wt w_1 \wt S w_2, \qquad 
	\wt w_1 \preceq \rho_{\varphi_1-\nu}(\iw_n) \quad \text{and} \quad \wt S = 0^{s+\nu}, 
\end{equation*}
and a~priori there are no reasons to say if $v$ coincide with $\wt v$ or not? 
If~we interpret the area covered by~$v$ as a window we monitor 
the sequence $\rho_{\om_n(j)}(\iw_n)$ through then it looks like 
the block $\rho_{\varphi_1}(\iw_n)$ is ``melting'' from the right side. 
We~erase $\nu$~letters in the subword $w_1$ while extending the spacer. 
Furthermore, from the statistical point of view 
the words $v$ and~$\tilde v$ are asymptotically $\dbar$-close with $\mu$-probability $1-\eps_n$, 
where ${\eps_n \to 0}$, if spacers are supposed to be asymptotically small. 
\medbreak

\proc{Definition.}
Let us consider the order on finite words: $w_1 \ge w_2$ if $w_1$ differs from $w_2$ by 
replacing some letters to the spacer symbol ``0''. We call {\it saturated complexity\/} 
of an infinite word~$w$ the function $\sC_w(l)$ counting the the total number of maximal points 
with respect to the order ``$\ge$'' in the set $\Lang(w,l)$ of subwords in $w$ of length~$l$. 
\medbreak

%

\begin{theo}
${\sC_{\iw_{n+1}}(h_n+1) \le h_n^3}$ for any iceberg system with spacers. 
\end{theo}

This theorem follows directly from the discussion of the privious remark.


\section{Lower bounds for symbolic complexity}
\label{sLowerBounds}

Now we pass to the discussion of lower estimates for $\p_{\iw_\infty}(l)$. 
The idea is to observe that generically, for iceberg systems with 
sufficiently rich combinatorial structure almost all configurations 
\begin{equation*}
	v = w_1 w_2 \preceq \rho_{\varphi_1}(\iw_n)\,\rho_{\varphi_2}(\iw_n), \qquad 
	w_1 \preceq \rho_{\varphi_1}(\iw_n), \qquad 
	w_2 \preceq \rho_{\varphi_2}(\iw_n), 
\end{equation*}
and, respectively, almost all triples $(|w_1|,\phi_1,\phi_2)$ 
are observed inside the word $\iw_{n+1}$. 
The following lemma is a~simple combinatorial observation. 

\begin{lemma}\label{lemOverlappingIntervals}
Let us consider two large 
intervals on the integer line $[k,k+m]$ and $[l,l+m]$ 
of the same length~$m$, and assume that $l-k \not= 0$. 
There exists a subset ${A \subset [k,k+m]}$ satisfying conditions 
\begin{equation*}
	A \cap (l+k-A) = \eset \quad \text{and} \quad \# A \ge \frac12\, m \,(1+o(1)). 
\end{equation*}
\end{lemma}

For the sequel it is convenient to find and fix a~value $m_0$ such that 
$\# A \ge \frac13\, m$ in the above lemma whenever ${m \ge m_0}$. 


Let us consider an iceberg system given by the independent random 
cyclic rotation parameters $\a_n(j)$ uniformely distributed on the integer interval $[0\ito h_n)$, 
where ${h_n = |\iw_n|}$. Let us use notation $u \cyclsub w$ in the case 
if $u$ is a~subword of some cyclic rotation of $w$, i.e.\ ${u \preceq \rho_\phi(w)}$ for some~$\phi$. 
This new class of dynamical systems generalizes the class or randomized rank one dynamical systems 
introduced by D.\,Ornstein (\cite{O}, see also \cite{AbPaPr}). 

\proc{Definition.}
We say that our symbolic system satisfies property $\cD_m(\beta)$ with ${\beta > 0}$ if 
any two subwords ${u \cyclsub \iw_n}$ and ${v \cyclsub \iw_n}$ of the same length ${m \ge \beta|\iw_n|}$ 
are always identical whenever they are equal. Here we use the word {\it identical\/} 
to define the case when two subwords has the same starting position in~$w$, i.e., in fact, 
they corresponds to one entrance of a~subword. For example, two entrances of the word ${u = v = cat}$ 
in a~lerger word ${w = littlecatandbigcat}$ starting at position $6$, and position $15$ 
are equal but not identical. 
\medbreak

In the sequel we will use symbol $\Prob(A)$ for the probability of the set $A$ according 
to the probability space hosting the random parameters $\om_n(j)$. 
It is convinient to consider two letter alphabet ${\Set{A} = \{0,1\}}$. 
Assume that the initial word in the iceberg construction $\iw_1$ contains $\nu\cdot|\iw_1|$ ones 
and $(1-\nu)\cdot|\iw_1|$. It follows easily from the definition that 
any word $\iw_n$ contains exactly $\nu$-fraction of ones. 
The following lemma mainly concerns the $n$'th step of the construction when 
we build the new word $\iw_{n+1}$ combining a~sequence of $q_n$ random cyclic rotations of~$\iw_n$. 

\begin{lemma}
\label{lemTrivialMatching}
Consider a pair of independent cyclic rotations, for example, 
\begin{equation*}
	W_1 = \rho_{\om_n(j)}(\iw_n) \quad \text{and} \quad W_2 = \rho_{\om_n(j+1)}(\iw_n) 
\end{equation*}
and a pair of intervals 
${I_1,I_2 \subset \{0,1,\dots,|\iw_n|-1\}}$. 
Let ${v_1 = W_1|_{I_1}}$ and ${v_2 = W_2|_{I_2}}$ be the two subwords corresponding 
to the subsets of indexes $I_1$ and $I_2$. Then the probability of matching ${v_1 = v_2}$ 
is less than $1 - 2\nu + 2\nu^2$. 
\end{lemma} 

\proc{Proof.}
It is enough to show that the conditional probability of matching $v_1 = v_2$ 
with respect to the finite set of random parameters that influenced to the structure of $\iw_n$. 
The probility of matching is less than the probability of matching of letters 
at only one posotion in $I_1$. The uniform random rotation generates at any position 
the random variable which is distributed 
as `$0$' with probability $1-\nu$ and `$1$' with probability~$\nu$. 
Thus, the one-letter matching process is equivalent to matching of two independent 
varibles given by the vector $\{1-\nu,\:\nu\}$. 
\ep
\medbreak

Remark that this lemma can be significantly strengthened, and our strategy in this paper 
is to avoid complicated investigation of random processes including large deviation mechanism 
and to use only elementary probability technique. 

\begin{lemma}
\label{lemDnChain}
Assume that $\beta > \max\{36,\:m_0+1\}\cdot q_n^{-1}$. Then 
\begin{equation*}
	\Prob\bigl(\cD_{n+1}(\beta) \mid \cD_n(1/2)\bigr) \ge 1 - q_n^3 \cdot h_n^{-\beta q_n/12}. 
\end{equation*}
\end{lemma}

\proc{Proof.}
Suppose that $\cD_m(1/2)$ is satisfied. 
Consider two subwords $v_1$ and $v_2$ of $\iw_{n+1}$ 
(or its cyclic rotation, further we will omit such kind of remark) starting at position
$k_1$ and $k_2$ respectively. Assume also that 
\begin{equation*}
	|v_1| = |v_2| = M \ge \beta|\iw_{n+1}| = \beta q_n h_n. 
\end{equation*}
Both $v_1$ and $v_2$ are of the form: 
\begin{equation*}
	v_i = v_i^{(head)}\,\rho_{\phi_{i,1}}(\iw_n)\cdots\rho_{\phi_{i,N_1}}(\iw_n)\,v_i^{(tail)}. 
\end{equation*}
In other words, it starts and ends with two incomplete blocks, enclosing the sequence 
of complete blocks $\rho_{\phi_{i,j}}(\iw_n)$. In order to understand if $v_1$ and $v_2$ are matched, 
${v_1 = v_2}$, let us write both $v_1$ and $v_2$ starting at zero position. 
At this point we must consider two cases which are essentially different. 
The first case is related to the pair of subwords which are very ``close'', that is 
the first positions $y_1$ and $y_2$ of $v_1$ and $v_2$ satisfy inequality 
\begin{equation*}
	|y_1 - y_2| \le \frac12 |\iw_n|. 
\end{equation*}
And the second case deals with the situation when matching $v_1$ and $v_2$ actually 
we match pairs of independently rotated blocks. 

In the first case we can simply apply axiom $\cD_n(1/2)$ to the situation when 
a~word $w = \rho_\phi(\iw_n)$ is matched with itsels shifted by non-zero amount of positions, 
since the size of overlapping is equal or grather than $|\iw_n|-|y_1-y_2| \ge |\iw_n|/2$. 


Concerning the second case
it~can be easily seen that $v_1$ and $v_2$ contain a sequence of overlapping block pairs 
\begin{equation*}
	(\rho_{\om_n(j)}(\iw_n), \rho_{\om_n(j+s)}(\iw_n)), \qquad j = j_0,\dots,j_0+m-1,  
\end{equation*}
satisfying the following conditions: 
\begin{equation*}
	\bigl| \rho_{\om_n(j)}(\iw_n) \cap \rho_{\om_n(j+s)}(\iw_n) \bigr| \ge \frac12 |\iw_n|, \qquad 
	m \ge \frac12 \frac{M}{|\iw_n|} \ge \frac12 \beta q_n.    
\end{equation*}
In this formulas we shortly write $j+s$ instead of $j+s \pmod{q_n}$, where 
$q_n$ is the total number of subblocks $\rho_{\om_n(j)}(\iw_n)$ in the rotated word $\iw_{n+1}$. 
Now let us apply lemma~\ref{lemOverlappingIntervals} to the intervals 
$[j_0,j_0+m-1]$ and $[j_0+s,j_0+s+m-1]$ on~$\Set{Z}$ and find a set $J$ 
such that ${\#J \ge 1/3 \cdot m}$ and $J \cap (s+J) = \eset$. 
Observe that the pair of random variables $(\om_n(j),\om_n(j+s))$ for ${j \in J}$ 
are globally independent, since all sets $\{j,j+s\}$, ${j \in J}$, are disjoint. 
Now applying axiom $\cD_n(1/2)$ we come to the following conclusion. 
The only possibility for a~pair of overlapping blocks 
$\rho_{\om_n(j)}(\iw_n)$ and $\rho_{\om_n(j+s)}(\iw_n)$ to match is to encounter 
the situation when these two blocks are proper rotated and the word in the overlapping 
enter identically in both blocks, and the probability of such events is exactlty $h_n^{-1}$. 
Thus, counting all the pairs of subwords $v_1$ and $v_2$ under testing, we have 
\begin{multline*}
	\Prob\bigl(\cD_{n+1}(\beta) \mid \cD_n(1/2)\bigr) \ge 1 - h_{n+1}^3 \cdot h_n^{-\#J} \ge \\ 
	\ge 1 - h_{n+1}^3 \cdot h_n^{-\beta q_n/6} \ge 1 - q_n^3 \cdot h_n^{-\beta q_n/12}. 
\end{multline*}
\ep
\medbreak

Suppose that some decreasing sequence $\beta_n \to 0$ is given. 
Applying the previous lemma we can easily find a~sequence $q_n$ such that 
for any $n$ the probability discussed in the lemma is positive and 
we can take for any step of the construction an appropriate configuration 
$(\om_n(0),\dots,\om_n(q_n-1))$ such that $\cD_n(\beta_n)$ is true for all~$n$. 

Though, remark that this reasoning still has a strong logical gap. 
In~fact, to start this process we need a~starting word $\iw_1$ satisfying $\cD_1(1/2)$ ``a~priori'', 
and actually we can do it by choosing such ``never matching'' word like 
\begin{equation*}
	10, \qquad 100, \qquad 1001, \quad \dots, 
\end{equation*}
but in view of futher applications to the isomorphism problem 
we need a version of this lemma that can be applied to {\em any\/} starting word~$\iw_1$.

\begin{lemma}\label{lemStrongMatchProbEst}
Suppose that $\beta > \max\{36,\:m_0+1\}\cdot q_n^{-1}$ and assume that the initial word 
of the construction contains $\nu$-fraction of symbol \/$1$. Then 
\begin{equation*}
	\Prob\bigl(\cD_n(1/2)\bigr) \ge 1 - 2q_n^3 \cdot \eta^{-\beta q_n/12}, 
	\label{eStrongMatchProbEst}
\end{equation*}
where $\eta = 1-2\nu+2\nu^2$. 
\end{lemma}

\proc{Proof.} 
Let us examine once more the first case in the proof of lemma~\ref{lemDnChain} 
taking into account that now we cannot say nothing about the matching of large blocks 
in~$\iw_n$. An example of such word poorly amenable to decoding is 
$
	1010101010\dots10. 
$

Without loss of generality we can assume that in a~pair of matched subwords $v_1$ and $v_2$ 
the second subwords starts just one position to the right of (the beggining of) $v_1$. 
In other words, if $y_2 = y_1 + 1$. The only iformation we have is just that 
the last symbol of any full block $\rho_{\om_n(j)}(\iw_n)$ in $v_2$ is compared 
to the first symbol of the next full block $\rho_{\om_n(j+1)}(\iw_n)$ included in $v_1$, 
and the probability of this event is less than $\eta$ 
(see lemma~\ref{lemTrivialMatching}). The additional summand to the probability of matching 
is estimated as
\begin{equation*}
	h_{n+1}^2 h_n \cdot \eta^{\#J} \le q_n^2 \cdot \eta^{-\beta q_n/12}. 
\end{equation*}
Further, if $|y_2-y_1| > |\iw_n|/2$ then the estimate is given by the value
\begin{equation*}
	h_{n+1}^3 h_n \cdot \eta^{\#J} \le q_n^3 \cdot \eta^{-\beta q_n/12} 
\end{equation*}
and the result follows. 
\ep
\medbreak

The idea of the forthcoming discussion is to examine, modulo axiom~$\cD_n(\beta_n)$, 
how many subwords $v$ of length $h_n$ can appear asymptotically in the word~$\iw_{n+1}$. 
Matching axioms $\cD_n$ imply that counting subwords is, in a sense, equivalent 
to counting configuration $(|w_1|,\phi_1,\phi_2)$ describing how $v$ covers 
two adjacent blocks $\rho_{\om_n(j)}(\iw_n)\,\rho_{\om_n(j+1)}(\iw_n) \preceq \iw_{n+1}$. 
Let us consider triples $(\xi,\phi_1,\phi_2)$ indexing such kind of configurations. 
Recall that for any subword $v$ in~$\iw_{n+1}$ of length $h_n$, 
\begin{equation*}
	v \preceq \rho_{\om_n(j)}(\iw_n)\,\rho_{\om_n(j+1)}(\iw_n). 
\end{equation*}
Here $\phi_1 = \om_n(j)$ and $\phi_2 = \om_n(j+1)$ are the corresponding cyclic rotation 
parameters of the adjacent blocks covering~$v$, and $\xi$ has the meaning of point 
separating these blocks inside~$v$. Notice that any tripple defines in a~unique way 
a~word $v$ that can occur as a~subword of~$\iw_{n+1}$. Let us denote this word as
\begin{equation*}
	v = V(\xi,\phi_1,\phi_2). 
\end{equation*}
Our purpose is to find a set of configurations such that the corresponding words 
$V(\xi,\phi_1,\phi_2)$ are different.

\begin{lemma}\label{lemTripleMatching}
Consider two tripes $(\xi,\phi_1,\phi_2)$ and $(\eta,\psi_1,\psi_2)$ and suppose that 
\begin{gather*}
	\beta_n h_n < \xi < \eta < (1-\beta_n)h_n, \\ 
	|\xi-\eta| > \beta_n h_n, 
\end{gather*}
and axioms $\cD_n(\beta_n)$ are satisfied. Then 
$V(\xi,\phi_1,\phi_2) = V(\eta,\psi_1,\psi_2)$ if and only if 
for any interval among $[0,\xi)$, $[\xi,\eta)$ and $[\eta,h_n)$ 
the corresponding subwords are matched identically, in other words, if 
\begin{equation*}
	\eta-\xi = \psi_1-\phi_1 = \psi_2-\phi_2 = h_n + \psi_1 - \phi_2. 
\end{equation*}
\end{lemma}

Observe that if the conditions of lemma~\ref{lemTripleMatching} are true 
then evidently $V(\xi,\phi_1,\phi_2) \not= V(\eta,\psi_1,\psi_2)$, except one case: 
${\phi_1 = \phi_2}$ and ${\psi_1 = \psi_2}$ which means that 
the adjacent blocks in the iceberg construction are cyclicly rotated 
by the same amount of positions. In this case the boundaries given by 
the positions $\xi$ and $\eta$ are transparent and cannot be recognized, since 
${\rho_\phi(ww) = \rho_\phi(w)\rho_\phi(w)}$. 

\proc{Proof of lemma~\ref{lemTripleMatching}.} The idea of this lemma is to require that 
all the intervals that appear while matching the words generated by 
these two configurations must be sufficiently long, so that we can apply 
property $\cD_n(\beta_n)$. In fact, for intervals $[0,\xi)$, $[\xi,\eta)$ and $[\eta,h_n)$ 
we know that the corresponding subwords in the blocks must be identically 
located subwords, and it is easy to represent this conclusion 
in terms of parameters $\xi$, $\eta$, $(\phi_1,\phi_2)$ and $(\psi_1,\psi_2))$. 
\ep
\medbreak

\begin{lemma}
Suppose that all the triples $(\xi,\phi_1,\phi_2)$ are found 
for subwords in $\iw_{n+1}$ of length~$h_n$. 
Then counting configurations with $\phi_1 \not= \phi_2$ we have 
\begin{equation*}
	\p_{\iw_{n+1}}(h_n) \gtrsim \frac{1-4\beta_n}{\beta_n} \cdot (h_n^2 - h_n) 
\end{equation*}
if $\beta_n \ge h_n^{-1+\delta}$ and ${0 < \delta < 1/2}$. 
\end{lemma}

\proc{Proof.}
Let us choose a progression $\kappa,\; \kappa+a,\; \dots,\; \kappa+(m-1)a$ 
of length $m$, where 
\begin{gather*}
	\kappa = a = \lfloor \beta_n h_n \rfloor, \\ 
	m = \left\lfloor \frac{(1-2\beta_n)h_n}{a} \right\rfloor \ge 
			\frac{(1-2\beta_n)h_n}{\beta_n h_n+1} - 1 \ge 
			\frac{(1-2\beta_n)(1-\beta_n^{(1-\delta)/\delta})}{\beta_n} -1 \ge \frac{1-4\beta_n}{\beta_n}, 
\end{gather*}
sinse $(1-\delta)/delta \ge 1$ for $\delta \le 1/2$. 
Any pair $(\tau_1,\tau_2)$ in the set of triples 
\begin{equation*}
	\Sigma = \bigl\{ (\kappa + ja,\phi_1,\phi_2) \where \phi_1 \not= \phi_2,\ 0 \le j < m \bigr\} 
\end{equation*}
satisfies the conditions of lemma~\ref{lemTripleMatching} and in addition we can state that 
${V(\tau_1) \not= V(\tau_2)}$, since ${\phi_1 \not= \phi_2}$. Finally, let us 
count the triples in~$\Sigma$: 
\begin{equation*}
	\#\Sigma = m \cdot (h_n^2-h_n) \ge \frac{1-4\beta_n}{\beta_n} \cdot (h_n^2 - h_n). 
\end{equation*}
\ep
\medbreak


\begin{theo}
Suppose that $\eps > 0$ is given, and $q_n$ satisfies the 
asymptotics $q_n \sim h_n^{\gamma}$ with $\gamma > 1 - \eps$. 
Then our iceberg system given by the sequence of i.i.d.\ random parameters 
$\om_n(j)$ with ${|\iw_{n+1}| = q_n\,|\iw_n|}$ almost surely generates 
an infinite word $\iw_\infty$ such that for ${n \ge n_0}$ 
\begin{equation*}
	\p_{\iw_\infty}(h_n) \ge \p_{\iw_{n+1}}(h_n) \ge h_n^{3-\eps} 
\end{equation*}
for ${n \ge n_0}$. 
\end{theo}

\proc{Proof.}
Let us take $\beta_n = h_n^{-1+\delta}$ with $\max\{1-\gamma,\,0\} < \delta < \eps$. 
In order to apply lemma~\ref{lemStrongMatchProbEst} and establish property $\cD_{n}(\beta_n)$ 
let us look at the probability in~\ref{eStrongMatchProbEst}: 
\begin{equation*}
	\Prob\bigl(\cD_n(1/2)\bigr) \ge 
		1 - 2q_n^3 \cdot \eta^{-\beta q_n/12} = 1 - \cE_n 
\end{equation*}
and
\begin{equation*}
	\cE_n = 2\eta^{3\gamma\log_\eta h_n - h_n^{(-1+\delta+\gamma)}/12}, 
\end{equation*}
where $\beta_n q_n = h_n^{-1+\delta+\gamma}$ and $-1+\delta+\gamma > 0$. 
%
%
Since the series $\sum_n \cE_n$ converges, we can apply Borel--Cantelli lemma 
and see that almost surely there exists $n_1$ such that $\cD_{n}(\beta_n)$ for ${n \ge n_1}$. 
The symbolic complexity is now estimated as follows:
\begin{equation*}
	\p_{\iw_\infty}(h_n) \ge \p_{\iw_{n+1}}(h_n) 
		\ge \const \cdot \frac{h_n^2}{\beta_n} 
		= \const \cdot h_n^{3-\delta} \ge h_n^{3-\eps}, 
\end{equation*}
whenever $n \ge n_0 \ge n_1$. 
\ep
\medbreak

This theorem can be strengthen as follows to get symbolic complexity 
arbitrary close to~$l^3$. 

\begin{theo}
Consider a sequence $q_n \sim h_n^{\gamma}$, $\gamma \ge 1$, 
defining a set of random iceberg systems given by the independent 
and uniformely distributed random parameters $\om_n(j)$. 
Then almost surely the symbolic complexity for this iceberg system satisfies the inequality 
\begin{equation*}
	\p_{\iw_\infty}(h_n) \ge h_n^{3-\eps}, \qquad n \ge n_0(\eps). 
\end{equation*}
\end{theo}

The only difference with the proof of the previous theorem is 
that we have to fix some maximal value $\eps_0$ and a~universal sequence 
$\beta_n = h_n^{-1+\delta_n}$ such that ${\delta_n \to 0}$ sufficiently slow. 

It was observed by S.\,Ferenczi \cite{Ferenczi2} that the lowest rate of growth 
for the complexity function conserning rank one systems can be seen along 
the subsequence $~h_n$. Thus, this case is concerned as the most difficult if we are interested in 
the estimate from below (see also lemma~\ref{lemUpperBoundIcebergSC}). 
Now, the above theorems can be easily extended to {\em all} lengths of subwords ${l \in \Set{N}}$. 
The complexity outside the sequence $h_n$ becomes even larger (cf.~\cite{Ferenczi2}), 
since we have more freedom combinating a word from several rotated copies of~$\iw_n$. 
Then we come to the following result. 

\begin{theo}\label{thCubicGrowthSC}
Given a random iceberg systems with $q_n \sim h_n^{\gamma}$, $\gamma \ge 1$, 
the infinite word $\iw_\infty$ generated by the system with probility one possesses 
almost cubic estimate for the complexity, 
\begin{equation*}
	\p_{\iw_\infty}(l) \ge l^{3-\eps}, \qquad n \ge n_0(\eps), 
\end{equation*}
where the function $n_0(\eps)$ depends on the system. 
\end{theo}

\section{Application to the calculation of rank}

Using the effects discussed in section~\ref{sLowerBounds} we can show that 
certain iceberg systems are not rank one. We give only an idea of proof. 

\begin{theo}
Keeping the conditions of theorem~\ref{thCubicGrowthSC} almost surely 
the iceberg system is not rank one. 
\end{theo}

Remark that the method we use is not applicable to the multiple rank property, for example, rank two. 

\proc{Sketch of the proof.}
Suppose that our transformation $T$ is rank one. 
Then there exists a partition $P = \{B_0,B_1\}$ 
and an orbit $x_k = T^k x_0$ such that the complexity of the word $(x_k)$ 
is quadratic along a sequence $H_i$ (see~\cite{Ferenczi2}),  
\begin{equation*}
	\p_{(x_k)}(H_i) \le \frac12 H_i^{2}. 
\end{equation*}
The most difficult point in the proof is to manage infinite words 
which are generated by {\em arbitrary\/} finite partition~$P$. 
In~other words, the following argument is used to show that the complexity 
of $(x_k)$ is $l^{3-\eps}$ providing a~contradiction. Having small $\zeta$, 
erasing $\zeta$-fraction of letters the following estimate remains true:  
\begin{equation*}
	\p_{\iw_{n+1}}(l) \ge l^{3-\eps}. 
\end{equation*}
Indeed, our infintie word $(x_k)$ is covered up to a small error $\zeta_n$ 
by cyclic rotations of $\iw_n$ and asymptotically it is easy to see from approximation
that passing through the iceberg with index $n$ names for all points (but $\zeta_n$-part) 
are $\zeta_n$-close. 

Applying these arguments we must take into account that $\zeta_n$ is small 
but it can go to zero arbitrary slowly. 
\ep
\medbreak

\section{Scaling approximation and Pascal adic transformation} 

Let us consider a tree associated with the Pascal triangle given by the vertice set 
$\cV = \{(k,n) \where 0 \le k \le n\}$, and the set of admissible paths in this graph: 
\begin{equation*}
	X = \bigl\{ (y_n,n) \where y_{n+1} = y_n\ \ \text{or}\ \ y_{n+1} = y_n+1 \bigr\}. 
\end{equation*}
Next, let us consider an order on $X$ induced by the natural order on the verticies on same level: 
$(k,n) \le (k+1,n)$. Pascal adic transformation $T$ is the transformation on $X$ 
that maps a path $x$ to the smallest path grater than $x$ (up~to a~small neglectible set). 
Pascal map $T$ becomes ergodic measure preserving transformation if we endow $X$ 
with the natural measure describing statistics of paths (see~\cite{VershikOnPascalContSp}). 

We mention Pascal transformation in this paper by two reasons. 
First, it belongs to the class of measure preserving maps containing iceberg systems as well 
and based on the concept of scaling approximation introduced below in this section. 
The second reason is the following theorem establishing the exact cubic asymptotics 
for the symbolic complexity of~$T$ (see~\cite{MelaPetersenOnDynPrPascal}). 

\begin{theo}(X.\,M\'ela, K.\,Petersen)
The symbolic complexity of the Pascal adic transformation $T$ satisfies $p(l) \sim \frac16\, l^3$. 
\end{theo}

\proc{Definition.}
We say that a measure preserving invertible map satisfies the property 
of {\it scaling approximation of rank one\/} with a scaling function $\scale(h)$ 
(or simply {\it scaling rank one\/}) if for any ${\eps > 0}$ 
almost every orbit encoded with a finite measurable partition 
can be $\eps$-covered by subwords $u_j$ of a fixed word $|W_{(\eps)}| = h$ 
such that the average length of $u_j$ is greater than $\scale(h)\,(1+o(1))$. 

The property of {\it funny scaling approximation of rank one}, by analogy with 
the funny rank one property, is definied in the same manner but taking instead of subwords 
arbitrary restrictions $v_j = W_{(\eps)}|_{I(j)}$ of the words $W_{(\eps)}$ 
considered as a function ${ W_{(\eps)} \Maps [0\ito |W_{(\eps)}|-1) \to \Set{A} }$. 
\medbreak

\begin{figure}[th]
  \centering
  \unitlength=1mm
  \includegraphics{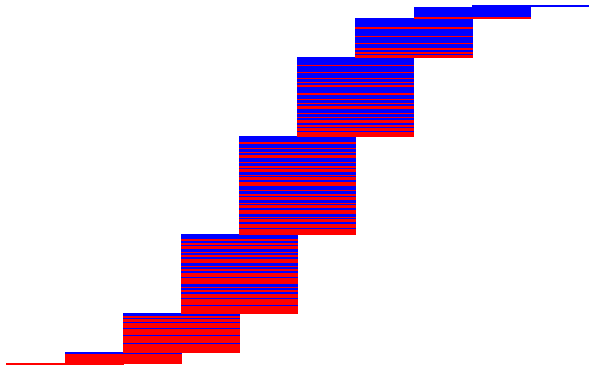} 
  \caption{On this figure we put the collection of towers in the cutting-and-stacking procedure 
  	for Pascal adic map, drawn in such a way that all the elementary sets sharing the same level 
  	are always labeled with the same letter. Here the $9$'th step od the construction is shown, 
  	the vertical columns are the Rokhlin towers of height $C(n,k)$, where ${n = 9}$ and ${k=0,1,\dots,n}$, 
  	and the common height of this generalized tower is equal to ${h_9 = 2^8 = 256}$. 
  	Red lines correspond to symbol `$1$' and blue lines -- to symbol `$1$'. 
  	$\;$ {\it Note\/}: you can zoom this~picture and see in details any part of the coding sequence 
  	$x_0,x_1,\dots,x_{255}$. }
  \label{fPascalScalingAppr}
\end{figure}

\begin{theo}
Any rank one (respectively rank~$m$) transformation has the property of scaling approximation 
of rank one (rank~$m$) with ${\scale(h) = 1}$. 
\end{theo}

\begin{theo}
Any iceberg map is of scaling rank one with ${\scale(h) = 1/2}$. 
\end{theo}

This theorem follows from the simple observation that any cyclic rotation 
of the word $\iw_n$ is a kind of interval exchange involving two parts, 
one of length $a|\iw_n|$ and another of length $(1-a)\iw_n$ which are 
compensated: $\frac12(a+(1-a)) = \frac12$. 

\proc{Observation 1.}
In can be easily seen that if $\scale(h) \ge \const > 0$ then the map has positive local rank 
(see~the scheme of the proof in~\cite{IcePaperI}). 
\medbreak

\proc{Observation 2.}
Pascal adic transformation admits scaling approximation of rank one with the scaling function
\begin{equation}
	\scale(h) = \frac1{\sqrt{\pi \,\log_2 h}}. 
	\label{ePascalApproximationScale}
\end{equation}
\medbreak

The picture on figure~\ref{fPascalScalingAppr} is the result of reconstruction of the usual 
cutting-and-stacking representation of the Pascal map, when the towers of height $C(n,k)$ 
are lifted so that the sets on the same level of the generalized tower 
are always marked by the same symbol. Here $C(n,k)$ are the binomial coefficients. 

\proc{Question 1.}
What is the precise symbolic complexity of the iceberg system with random rotations? 
\medbreak

\proc{Question 2 (iceberg systems with multiple IET).}
Let us consider a class of systems given by the more general construction: 
instead of cyclic rotation (which is in fact a kind of discrete two interval exchange map) 
we apply a~random $r$-interval exchange transformation to the words $\iw_n$ (see~\cite{IcePaperI}). 
What is the typical asymptotics for the symbolic complexity for this ensemble of symbolic systems?
\medbreak

\proc{Question 3.}
Is it true that the scaling function given by~\eqref{ePascalApproximationScale} 
is optimal for the Pascal adic transformation? 
\medbreak


\acks
The author is very greatful to M.\,Lemanczyk and S.\,Ferenczi for the discussions 
concerning the subject of this work, largely inspired by papers \cite{Ferenczi2} and \cite{FerencziLem}. 
 
I would like to thank A.\,St\"epin, A.\,Vershik, B.\,Gurevich, V.\,Oseledec, S.\,Pirogov, 
J.-P.\,Thouvenot, K.\,Petersen, V.\,Ryzhikov, El H.~El~Abdalaoui  
and the participants of seminar ``Ergodic theory and statistical physics'' at Moscow State University 
and ``St.~Petersburg seminar on representation theory and dynamical systems'' 
for the fruitful discussions and attention to this work. 

\bibliographystyle{plain}
\bibliography{IcePaperI}

\begin{thebibliography}{10}

\bibitem{AbPaPr}
El~H.~El Abdalaoui, F.~Parreau, and A.A. Prikhod'ko.
\newblock A new class of ornstein transformations with singular spectrum.
\newblock {\em Annales de l'Institut Henri Poincare (B) Probability and
  Statistics}, 42(6):671--681, 2006.

\bibitem{ArnouxMauduitShTam}
P.~Arnoux, C.~Mauduit, I.~Shiokawa, and J.-I. Tamura.
\newblock Complexity of sequences defined by billiards in the cube.
\newblock {\em Bull. Soc. Math. France}, 122:1--12, 1994.

\bibitem{BRL}
S.~Brlek.
\newblock Enumeration of factors in the thue-morse word.
\newblock {\em Discr. Appl. Math.}, 24:83--96, 1989.

\bibitem{delVAR}
A.~de~Luca and S.~Varrichio.
\newblock Some combinatorial properties of the thue-morse sequence.
\newblock {\em Theor. Comput. Sci.}, 63:333--348, 1989.

\bibitem{Ferenczi2}
S.~Ferenczi.
\newblock Rank and symbolic complexity.
\newblock {\em Ergod.\ Th.\ and Dynam.\ Sys.}, 16:663--682, 1996.

\bibitem{FerencziOnMTCompl}
S.~Ferenczi.
\newblock Measure-theoretic complexity of ergodic systems.
\newblock {\em Israel J. Math.}, 100:189--207, 1997.

\bibitem{Ferenczi1}
S.~Ferenczi.
\newblock Systems of finite rank.
\newblock {\em Colloq.\ Math.}, 73(1):35--65, 1997.

\bibitem{FerencziLem}
S.~Ferenczi and M.~Lemanczyk.
\newblock Rank is not a spectral invariant.
\newblock {\em Studia Math.}, 98(3):227--230, 1991.

\bibitem{FerencziZamboniOnIET08}
S.~Ferenczi and L.Q. Zamboni.
\newblock Languages of $k$-interval exchange transformations.
\newblock {\em Bull. Lond. Math. Soc.}, 40(4):705--714, 2008.

\bibitem{MelaPetersenOnDynPrPascal}
X.~Mela and K.~Petersen.
\newblock Dynamical properties of the pascal adic transformation, {\it arxiv:
  0310317}.
\newblock 2003.

\bibitem{O}
D.S. Ornstein.
\newblock On the root problem in ergodic theory.
\newblock {\em Proc.\ 6th Berkley Sympos.\ Math.\ Statist.\ Probab., Univ.\
  Calif.}, 2:347--356, 1970.

\bibitem{IcePaperI}
A.A. Prikhod'ko.
\newblock On ergodic properties of ``iceberg'' transformations.
  i:~approximation and spectral multiplicity, {\it preprint},
  arxiv:1002.2808v1.

\bibitem{VershikOnPascalContSp}
A.M. Vershik.
\newblock The pascal automorphism has a continuous spectrum.
\newblock {\em Funct. Anal. Appl.}, 45(3):173--186, 2011.

\bibitem{VershikOnScalingEntrAndPPS}
A.M. Vershik.
\newblock Scaling entropy and automorphisms with purely point spectrum.
\newblock {\em Algebra i Analiz}, 23(1):111--135, 2011.

\end{thebibliography}

\end{document}